\documentclass{amsart}
\usepackage{amsmath}
\usepackage{amsfonts}
\usepackage{amssymb}
\usepackage{graphicx}%
\newtheorem{theorem}{Theorem}[section]
\newtheorem{lemma}[theorem]{Lemma}
\newtheorem{proposition}[theorem]{Proposition}
\theoremstyle{definition}

\newtheorem{ass}{\textbf{Assumption}}

\theoremstyle{remark}
\newtheorem{remark}[theorem]{Remark}
\numberwithin{equation}{section}

 \allowdisplaybreaks[4]
\begin{document}

\title{A Donsker-type  Theorem  for Log-likelihood Processes}


\author{Zhonggen Su}

\address{ School of Mathematical Science, Zhejiang University, Hangzhou, China. }
\email{suzhonggen@zju.edu.cn}


\author{Hanchao Wang}
\address{Zhongtai Securities Institute for Financial Studies,  Shandong University,  Jinan,  250100, China.}
\email{wanghanchao@sdu.edu.cn}


\date{}


\begin{abstract}
Let $(\Omega, \mathcal{F}, (\mathcal{F})_{t\ge 0}, P)$ be a complete stochastic basis, $X$ a semimartingale with predictable compensator $(B, C, \nu)$. Consider a family of probability measures $\mathbf{P}=( {P}^{n, \psi}, \psi\in \Psi, n\ge 1)$, where $\Psi$ is an index set, $ {P}^{n, \psi}\stackrel {loc} \ll{P}$, and denote the  likelihood ratio process by $Z_t^{n, \psi} =\frac{dP^{n, \psi}|_{\mathcal{F}_t}}{d P|_{\mathcal{F}_t}}$.  Under some regularity conditions  in terms of logarithm entropy and Hellinger processes, we prove that $\log Z_t^{n}$ converges weakly to a Gaussian process in $\ell^\infty(\Psi)$ as $n\rightarrow\infty$  for each fixed $t>0$.
  \end{abstract}

\maketitle

\section{Introduction and  Preliminaries}

The celebrated Donsker theorem is a functional extension of the central limit theorem in probability theory. Plenty of research on this topic has come out in the past decades. The reader is referred to classic books and papers like Dudley \cite{d},  Gine and Zinn \cite{gz}, Ossiander  \cite{o}, Andersen et. al \cite{agoz}, Liptser and Shiryaev  \cite{ls}, van der Geer \cite{va}, Billingsley \cite{b},  Jacod and Shiryaev \cite{js} for both theoretical framework and wide applications. A primary purpose of the present paper is to establish a certain Donsker theorem for log-likelihood processes indexed by  an arbitrary set. In this section, we first  introduce some basic notions about log-likelihood processes and  martingale representation property.

Throughout this  paper, we follow the standard definitions and notations of martingale theory, which can be found in the book by Jacod and Shiryaev \cite{js}.  Let $(\Omega, \mathcal{F}, (\mathcal{F})_{t\ge 0}, P)$ be a complete stochastic basis. Fix  a semimartingale $X$ on it, and assume that all $P$-martingales have a representation property relative to $X$. Denote by the triplet $(B, C, \nu)$ the predictable characteristic of $X$ (associated to some bounded truncation function). More precisely, if $\Delta X_t=X_t-X_{t-}$ denotes the jump of $X$ at time $t$, then $X_t-\sum_{s\le t}(\triangle X_s-h_\tau(\triangle X_s))$ , where $h_{\tau}(x)=x1_{(|x|\le \tau)}$, is a special semimartingale, which can be uniquely divided into a bounded variation process and  a local martingale process. The $B$ is a bounded variation process of  $X-\sum_{s\le \cdot}(\triangle X_s-h_\tau(\triangle X_s))$.   Let $X^c$ be the continuous local martingale part of  $X$, then
\begin{eqnarray}
C_t=\langle X^c, X^c \rangle_t. \nonumber
\end{eqnarray}
Let $\mu$ be the jump measure of $X$ defined by
\begin{eqnarray}
\mu(\omega, dt, dx)= \sum_{s} \mathbf{1}_{(\Delta X_s(\omega)\neq 0)}\varepsilon_{(s, \Delta X_s(\omega))}(dt, dx)
\end{eqnarray}
where $\varepsilon_{(s, \Delta X_s(\omega))}$ denotes the Dirac measure at point $(s, \Delta X_s(\omega))$. The  $\nu$  is the unique predictable compensator of $\mu$ (up to a $P$-null set). Namely, $\nu$ is  a  predictable random measure such that for any predictable function\footnote{\,  Let $\tilde{\Omega}=\Omega \times \mathbb{R}_+\times \mathbb{R}$, $\tilde{\mathcal{P}}=\mathcal{P}\otimes \mathcal{B}$, where $\mathcal{B}$ is a Borel $\sigma$-field on $\mathbb{R}$ and $\mathcal{P}$ a $\sigma$-field  generated by all left continuous adapted processes  on $\Omega \times \mathbb{R}_+$. The predictable function is a $\tilde{\mathcal{P}}$-measurable function on $\tilde{\Omega}$.} $W$, $W\ast (\mu-\nu)=W\ast  \mu-W\ast\nu)$ is a  local martingale, where the $W\ast  \mu$ is defined by
\begin{eqnarray}
& &W\ast \mu_t(\omega)\nonumber\\
 & &= \left\{
\begin{array}{ll}
 \int_0^t\int_{\mathbb{R}}W(\omega, s,x)\mu(\omega; ds, dx)& \mbox{if} \int_0^t\int_{\mathbb{R}}|W(\omega, s,x)|\mu(\omega; ds, dx)<\infty, \\
 \\
 +\infty & \mbox{otherwise}.
 \end{array}
 \right. \nonumber\label{StarIn}
\end{eqnarray}
(See Section 2.1 of Chapter 2 in Jacod and Shiryaev \cite{js} for more details).   Note the predictable quadratic variation is given by
 \begin{eqnarray}
 \langle  W\ast (\mu-\nu)  \rangle_t = \big(W-\hat{W}\big)^2\ast \nu_t +\sum_{s\le t}(1-a_s)\hat{W}^2_s,
\end{eqnarray}
where
\begin{eqnarray}
  \hat{W}=\hat{W}(\omega, t) = \int_{\mathbb{R}} W(\omega, t, x)\nu(\omega, \{t\}\times dx)
\end{eqnarray}
and
\begin{eqnarray}
   a_t(\omega)=\nu(\omega, \{t\}\times \mathbb{R}).
\end{eqnarray}
It follows from Corollary 1.19 of Chapter 2 in Jacod and Shiryaev \cite{js} that $a_t=0$ is equivalent to the fact $X$ is a quasi-left continuous process\footnote{\,  A quasi-left continuous process $X$ is a   c\`adl\'ag adapted process such that for any   increasing stoping times $T_n$ with limit $T$, $\lim_{n\to \infty } X_{T_n}=X_T.$}.  Specially, for a  process  with independent increments, $a_t=0$ means this process  has no fixed time of  discontinuity \footnote{\, $t$ is called as the fixed time of discontinuity  if $P(\Delta X_t\neq 0)>0$.}. Thus we may and do choose a good version of both $\hat{W}$ and $a$ such that $\hat{W}$ is the predictable projection of $W(\omega, t, \Delta X_t)\mathbf{1}_{(\Delta X_t\neq 0)}$ and  $a_t\le 1$.  In particular,
  \begin{eqnarray}
    E\big(\mathbf{1}_{(\Delta X_t=0)}\big| \mathcal{F}_{t-}\big)=1-a_t.
\end{eqnarray}

 Now consider another probability measure $P'$  such that
\begin{eqnarray}
   P'\stackrel {loc} \ll P,
\end{eqnarray}
which means that for any $t\ge 0$, $P'|_{\mathcal{F}_t}\ll P|_{\mathcal{F}_t}$. Define the likelihood ratio process
\begin{eqnarray}
   Z_t=\frac{dP'|_{\mathcal{F}_t}}{dP|_{\mathcal{ F}_t}}.
\end{eqnarray}
It follows from Chapter III  in Jacod and Shiryaev \cite{js} that $Z_t$ is a local martingale.

Since  by assumption all $P$-martingales have a representation property relative to $X$, according to Theorem 5.19 of Chapter III in Jacod and Shiryayev \cite{js}, $Z_t$ has the following representation: there is a predictable process $\beta$ and a nonnegative predictable function $Y$ on $\tilde{\Omega}$ such that
\begin{eqnarray}\label{nn1}
   Z_t=\left\{
   \begin{array}{ll}
   Z_0 e^{N_t-\frac 12 \beta^2\cdot C_t} \prod_{s\le t}(1+\Delta N_s)e^{-\Delta N_s}, &\quad (\omega, t)\in \triangle, \\
   0, &\quad \mbox{otherwise}.
   \end{array}
   \right.
\end{eqnarray}
Here
\begin{eqnarray*}
   N_t=\beta\cdot X_t^c +\Big(Y-1+\frac{\hat{Y}-a}{1-a}\mathbf{1}_{(a<1)}\Big)\ast (\mu-\nu)_t
\end{eqnarray*}
and $\triangle$ is a random set defined as follows
\begin{eqnarray*}
   \sigma=\inf\{t: \hat{Y}_t>1 \mbox{ or } \hat{Y}_t<1 \mbox{ and } a_t=1\},
\end{eqnarray*}
\begin{eqnarray*}
    H_t =\frac{1}{2}\beta^2\cdot C_t + \big(1-\sqrt Y\big)^2\mathbf{1}_{[0, \sigma)}\ast \nu_t+ \sum_{s\le t}\Big(\sqrt{1-a_s}-\sqrt{1-\hat{Y}_s}\Big)^2\mathbf{1}_{(s<\sigma) },
\end{eqnarray*}
\begin{eqnarray*}
     T_n=\inf\{t: H_t\ge n\},
\end{eqnarray*}
\begin{eqnarray*}
     \triangle =[0,\, \sigma) \cap (\cup_{n}[0,\, T_n]).
\end{eqnarray*}
Note that $\beta, Y$ and $\triangle$ depend on $P'$.  In fact, $Y$ can be explicitly represented as follows. Let $M_{\mu}^{P}$ be a measure on $(\tilde{\Omega}, \tilde{\mathcal{P}})$ where $\tilde{\Omega}:=\Omega \times \mathbb{R}_+\times \mathbb{R}$, $\tilde{\mathcal{P}}:=\mathcal{P}\otimes \mathcal{B}$, such that $M_{\mu}^{P}(W)=E(W*\mu)_{\infty}$ for all measurable nonnegative functions $W$. Then  $Y$ is  the conditional expectation of $\frac {Z}{Z_{-}}$ with respect to $\tilde{\mathcal{P}}$ under $M_{\mu}^{P}$, namely
$$Y =M_{\mu}^{P}\Big(\frac {Z}{Z_{-}}|\tilde{\mathcal{P}}\Big).$$

Define the log-likelihood process $L_t$ by
\begin{eqnarray}
      L_t = \log Z_t.
\end{eqnarray}
This process has been a well-studied object in the context of  both stochastic processes and statistical inferences.  Obviously,
\begin{eqnarray*}
      L_t &=&  \log Z_0+\beta\cdot X_t^c -\frac{1}{2} \beta^2\cdot C_t + \Big(Y-1+ \frac {\hat{Y}-a}{1-a}\mathbf{1}_{(a<1)}\Big)\ast(\mu-\nu)_t\\
       && + \sum_{s\le t}\big(\log (1+\Delta N_s)-\Delta N_s\big)\nonumber
\end{eqnarray*}

Assume we are given a family of probability measures $\mathbf{P}^n=\{P^{n, \psi}: \psi\in \Psi\}$ on $(\Omega,\mathcal{F})$, indexed by an arbitrary non-empty set $\Psi$, and assume
 $$P^{n,\psi}\stackrel{\text{loc}}{\ll}P $$
 for every $n>0$ and $\psi\in \Psi$. We shall be mainly interested in the sequence of  likelihood ratio processes $Z_t^{n, \psi}$. The main purpose of the paper is to establish a certain  Donsker theorem for log-likelihood processes $\log Z_t^n$ in $\ell^{\infty}(\Psi)$ as $n\rightarrow\infty$, where we denote by $\ell^{\infty}(\Psi)$  the space of bounded real-valued functions defined on $\Psi$.

 It seems hard to develop directly an invariance principle for $\log Z^n$ due to complicated structure. To the best of our knowledge, there are only few works in this area, such as Le Cam \cite{l}, Vostrikova \cite{v} and so on. The reader can find some interesting results in Nishiyama \cite{n1} and \cite{n3}, where $\log Z^n$ is assumed to be very special  continuous semimartingales and discrete time semimartingales respectively. It is a challenging problem to extend Nishiyama's work to general setting. To attack such a problem, we shall combine stochastic calculus techniques and chaining arguments with the Kakutani-Hellinger distance for probability measures. In particular, we shall  characterize the regularity  of $\ell^{\infty}(\Psi)$-valued log-likelihood processes in terms of the Kakutani-Hellinger distance  and the Hellinger processes.

The rest of the paper is organized as follows. We will first  make some necessary assumptions and then state our main results in the Section 2. The proof of main results is given in Section 3, which consists of several lemmas and two propositions.

\section{Main Results}
To state our main results, we need some more notations and make some technical assumptions. Start with the Kakutani-Hellinger distance between two probability measures $P$  and $P'$. Assume that $Q$ is a third probability measure on $(\Omega,\mathcal{F})$ such that
\begin{eqnarray}
 P \ll Q, \quad   P' \ll Q
\end{eqnarray}
Let
\begin{eqnarray}
     Z= \frac{dP}{dQ}, \quad   Z'= \frac{dP'}{dQ}
\end{eqnarray}
and define the Kakutani-Hellinger distance by
\begin{eqnarray}
     \rho^2(P, P') =\frac{1}{2}\int_{\Omega} \big(\sqrt Z- \sqrt{Z'}\big)^2 dQ.
\end{eqnarray}
 It is easy to check that $\rho (P, P')$ is a metric in the space of probability measures and does not depend on the probability measure $Q$.  Note
\begin{eqnarray}
     \rho^2(P, P') = 1- E_Q\sqrt{ZZ'}
\end{eqnarray}
For $0<\alpha<1$, call $\breve{H}(\alpha; P, P')= E_Q(Z^\alpha (Z')^{1-\alpha})$ the Hellinger integral of order $\alpha$. We remark that  $ \breve{H}(\alpha; P, P')\rightarrow 1$ as $\alpha\rightarrow 0$ if $P'\ll P$.

Proceed to introduce the Hellinger processes. Assume that
\begin{eqnarray}
      P\stackrel {loc}\ll Q, \quad P' \stackrel {loc}\ll Q
\end{eqnarray}
and define
\begin{eqnarray}
     Z_t=\frac{dP|_{\mathcal{F}_t}}{dQ|_{\mathcal{F}_t}}, \quad  Z'_t=\frac{dP'|_{\mathcal{F}_t}}{dQ|_{\mathcal{F}_t}}.
\end{eqnarray}
Then for each $0<\alpha<1$,  there is a unique predictable increasing process $h(\alpha; P, P')$, called the Hellinger process of order $\alpha$, such that

\noindent
(i)
\begin{eqnarray}
h(\alpha; P, P')_0=0\nonumber
\end{eqnarray}
  (ii)
 \begin{eqnarray}
 h(\alpha; P, P')_t= \mathbf{1}_{\cup[0,\, S_n]}\cdot h(\alpha; P, P')_t\nonumber
 \end{eqnarray}
  (iii)
 \begin{eqnarray}
  Y(\alpha)_t+Y(\alpha)_-\ast \nu_t \quad \mbox{ is local martingale }\nonumber
\end{eqnarray}
where
\begin{eqnarray*}
S_n=\inf\{t: Z_t>n \, \mbox{ or }\, Z_t'>n\}
\end{eqnarray*}
and
 \begin{eqnarray*}
 Y(\alpha)_t=Z_t^\alpha(Z_t')^{1-\alpha}.
\end{eqnarray*}

One can  extend the above Hellinger process to order zero and even to a general function. Given a function $\psi:\mathbb{R}\mapsto \mathbb{R} $ such that
\begin{eqnarray}
      \frac{\psi(x)}{|x-1|^2\wedge |x-1|}
\end{eqnarray}
 is bounded with convention $\frac{0}{0}=0$ and $\psi (1)=0$.  Denote
\begin{eqnarray}
     \jmath (\psi; P, P')_t =\sum_{s\le t}\frac{Z_s'}{Z'_{s-}}\psi\Big(\frac{Z_s/Z_{s-}}{ Z'_s/Z'_{s-}} \Big)
\end{eqnarray}
then there is a predictable increasing process, denoted by $\imath(\psi;  P, P')$, such that

\noindent
(i')
\begin{eqnarray*}
 \imath(\psi; P, P')_0=0 \nonumber
 \end{eqnarray*}
  (ii')
\begin{eqnarray*}
{\imath(\psi; P,P')}_t= \mathbf{1}_{\cup[0,\, T_n]}\cdot \imath(\psi; P, P')_t\nonumber
\end{eqnarray*}
(iii')
\begin{eqnarray*}
\jmath(\psi; P, P')_{T_n\wedge t}-\imath(\psi; P,P')_{T_n\wedge t} \quad \mbox{ is local martingale }.\nonumber
\end{eqnarray*}
Call $\imath(\psi; P,P')$ the Hellinger process of order 0 associated to $\psi$. In particular, if
\begin{eqnarray*}
     \psi(x)=\left\{
     \begin{array}{ll}
     1,& x=0,\\
     0, & x>0,
     \end{array}
     \right.
\end{eqnarray*}
then we simply call $\imath(\psi; P, P')$ the Hellinger process of order 0.

In general, it is rather complicated to compute $h(\alpha; P, P')$. However,  we fortunately have the following explicit formula in the special case  $P'\ll P$:

\begin{eqnarray*}
     h(\alpha; P, P')=\frac{\alpha(1-\alpha)}{2} \beta^2\cdot C+\varphi_\alpha(1, Y)\ast \nu + \sum_{s\le t} \varphi_\alpha(1-a_s, 1-\hat{Y}_s)
\end{eqnarray*}
In particular,
\begin{eqnarray*}
     h\Big(\frac 12; P, P'\Big)=\frac{1}{8} \beta^2\cdot C+\frac 12 (1-\sqrt Y)^2\ast \nu + \frac 12 \sum_{s\le t} \big( \sqrt{1-a_s}- \sqrt{ 1-\hat{Y}_s}\big)^2,
\end{eqnarray*}

\begin{eqnarray*}
     \imath(\psi; P, P')= Y \psi \Big(\frac 1Y\Big) \ast \nu + \sum_{s\le t} (1-a_s') \psi\Big(\frac{1-a_s}{1-a_s'}\Big).
\end{eqnarray*}

Our technical  assumptions mainly involve three aspects: the predictable envelope  of $\{Y^{n, \psi}, \psi\in \Psi\}$,  the Kakutani-Hellinger distance between probability measures $P^{n, \psi}$, and the size of index set $\Psi$.

For every $n>0$, denote the essence supremum $ \overline{Y}^{n}(\Psi)=[\sup_{\psi\in \Psi}Y^{n,\psi}]_{\tilde{\mathcal{P}}, M_{\nu}^P}$. This is  the predictable envelope of  $\{Y^{n,\psi}, \psi\in\Psi\}$  used in Definitions 2.1 and 2.3 of Nishiyama \cite{n3}.
\begin{ass}\label{a1}
 \it{For any $n>0$, $\psi\in \Psi$, $\triangle^{n,\psi}\equiv\Omega\times [0,1]$ and $0\le a<1$. Moreover, $\{Y^{n,\psi}, \psi\in \Psi\}$  attains their predictable envelope for every $n>0$, namely, there is a $\psi_0\in \Psi$ such that $Y^{n,\psi_0}=[\sup_{\psi\in \Psi}Y^{n,\psi}]_{\tilde{\mathcal{P}}, M_{\nu}^P}$.}
 \end{ass}
\begin{ass}\label{a2} \it{  For every  $\varepsilon>0$, as $n\rightarrow\infty $
\begin{eqnarray}
    \sup_{\psi\in \Psi}\imath \big(f_{1+\varepsilon},P,P^{n,\psi}\big)_{t}\stackrel{P }{\longrightarrow}0
     \end{eqnarray}
where $f_{1+\varepsilon}(x)=|x-1|1_{\{1/(1+\varepsilon)<x<1+\varepsilon\}^{c}}$.

\noindent
There is a nonnegative definite continuous function $V_t$ on $\Psi\times \Psi$, such that as $n\rightarrow\infty$,
 \begin{eqnarray}
 \sup_{\psi\in \Psi}\Big|h\big(\frac 12; P, P^{n,\psi}\big)_{t}-\frac{1}{8}V^{\psi,\psi}_t\Big|\stackrel{P }{\longrightarrow} 0
  \end{eqnarray}
  and  for every $\psi, \phi \in \Psi$,
   \begin{eqnarray}
   h\big(\frac 12; P^{n,\psi},P^{n,\phi}\big)_{t}\stackrel{P}{\longrightarrow} \frac{1}{8}V^{\psi,\phi}_t.
   \end{eqnarray}
}

\end{ass}
Let $\Psi$ be an arbitrary set,  $\Delta_\Pi$  a positive rational number. $\Pi=\{\Pi(\varepsilon)\}_{\varepsilon\in (0,\Delta_\Pi ]}$,  is called a decreasing
series of finite partitions (DFP) of $\Psi$ if

\noindent
          (i)  each $\Pi(\varepsilon)=\{\Psi(\varepsilon;k):1\le k\le N_{\Pi}(\varepsilon)\}$ is a finite partition of $\Psi$, namely
                      $$\Psi=\bigcup_{k=1}^{N_{\Pi}(\varepsilon)}\Psi(\varepsilon;k);$$

\noindent
          (ii) $N_{\Pi}(\Delta_{\Pi})=1$ and $\lim_{\varepsilon\downarrow 0}N_{\Pi}(\varepsilon)=\infty$;

\noindent
          (iii) $N_{\Pi}(\varepsilon)\ge N_{\Pi}(\varepsilon')$ as $\varepsilon\le \varepsilon'$.

 Given a $0<\varepsilon\le \Delta_\Pi$, the $\varepsilon$-entropy $H_{\Pi}(\varepsilon)$ is defined by
                         $$H_{\Pi}(\varepsilon)=\sqrt{\log (1+N_{\Pi}(\varepsilon))}.$$

\begin{ass}\label{a3}  \it{There exists a   decreasing
series of finite partitions, $\Pi$, of $\Psi$ such that as $n\rightarrow\infty$
$$\int_{0}^{\Delta_{\Pi}}H_{\Pi}(\varepsilon)d\varepsilon<\infty$$
 and
$$\|\mathfrak{H}^{n}\|_{\Pi }=O_{P}(1)$$
where
$$\|\mathfrak{H}^{n}\|^2_{\Pi}=\sup_{\varepsilon\in (0,\Delta_\Pi]\cap \mathbb{Q}}
\max_{1\le k\le
N_{\Pi}(\varepsilon)}\max_{\psi,\phi\in\Psi(\varepsilon,k)}\frac{1}{\varepsilon^2}
h\big(\frac 12; P^{n,\psi},P^{n,\phi}\big)_{1}\ .$$
}
\end{ass}
We are now ready to state our main result as follows.
\begin{theorem}\label{m1}
Under  Assumptions \ref{a1}, \ref{a2} and \ref{a3}, we have
    \begin{eqnarray}
    L^{n}_1\Rightarrow G \quad \mbox{ in } \ell^{\infty}(\Psi),  \label{LNG}
    \end{eqnarray}
 where $G$ stands for a Gaussian element in $\ell^{\infty}(\Psi)$, each $d$-dimensional marginal $(G^{\psi_{1}},\cdots,G^{\psi_{d}})$ is a  normal random vector with
mean
\begin{eqnarray}
\vec{\mu}=-\frac 12\big(V_1^{\psi_{i},\psi_{i}}, \, 1\le i\le d\big)
 \end{eqnarray}
and covariance structure
\begin{eqnarray}
 \Sigma=\big(V_1^{\psi_{i},\psi_{j}}\big)_{1\le i,j\le d}.
 \end{eqnarray}
\end{theorem}
The proof will be given in Section 2. For sake of comparison, we review an earlier  result due to Nishiyama \cite{n3} in the discrete time case. Let $  (\mathcal{F}_{i})_{i\ge 0} $ be a discrete time stochastic basis, and  $\mathbf{P}^{n}=\{P^{n,\psi}:\psi\in \Psi\}$  a family of probability measures on $(\Omega,\mathcal{F})$, such that
 \begin{eqnarray}
P^{n,\psi}\stackrel{\text{loc}}{\ll}P.
 \end{eqnarray}
Define
 \begin{eqnarray}
 W_{i}^{n,\psi}=\frac{d P^{n,\psi} |_{\mathcal{F}_i}}{d P|_{\mathcal{F}_i}}
  \end{eqnarray}
and
     \begin{eqnarray}
     \xi^{n,\psi}_{i}=\sqrt{\frac{W^{n,\psi}_{i}}{W^{n,\psi}_{i-1}}}-1
      \end{eqnarray}
Nishiyama \cite{n3} studied weak convergence for log-likelihood processes $\log W^n_n$ in $\ell^\infty(\Psi)$ and obtained a similar  result to (\ref{LNG}) under some integrability assumptions involving $\xi^n$'s and entropy conditions. More specifically, assume

\noindent
(i)  for every $\varepsilon>0$
\begin{equation*}
\sum_{i=1}^{n}E\Big( (\sup_{\psi\in  \Psi}\xi^{n,\psi}_{i})^{2}1_{\{\sup_{\psi\in \Psi}\xi^{n,\psi}_{i}>\varepsilon\}}\Big|\mathcal{F}_{i-1}\Big)\stackrel{P}{\longrightarrow}0;
\end{equation*}
(ii)  there exists a   decreasing
series of finite partitions, $\Pi$, of $\Psi$  such that
 \begin{eqnarray}
 \sup_{\varepsilon\in (0,\Delta_\Pi]\cap \mathbb{Q}}\max_{1\le k\le N_{\Pi}(\varepsilon)}\frac{1}{\varepsilon^2} \sum_{i=1}^{n}\sup_{\psi,\phi\in\Psi(\varepsilon,k)}E\big(|\xi^{n,\psi}_{i}-\xi^{n,\phi}_{i}|^{2}|\mathcal{F}_{i-1}\big) =O_{P}(1)
  \end{eqnarray}
and
 \begin{eqnarray}
 \int_{0}^{\Delta_\Pi }H_{\Pi}(\varepsilon)d\varepsilon<\infty;
 \end{eqnarray}

\noindent
(iii) there is a $V:\Psi\times \Psi\rightarrow R$ such that
  \begin{equation*}
  \sup_{\psi\in \Psi}\Big|\sum_{i=1}^{n}4E^{*}\big((\xi^{n,\psi}_{i})^{2}|\mathcal{F}_{i-1}\big)-V^{\psi,\psi}\Big|\stackrel{P}{\longrightarrow}0;
  \end{equation*}
and for   $\psi,\phi \in \Psi$
   \begin{equation*}
   \sum_{i=1}^{n}4E\big(\xi^{n,\psi}_{i}\xi^{n,\phi}_{i}|\mathcal{F}_{i-1}\big)\stackrel{P}{\longrightarrow} V^{\psi,\phi}.
  \end{equation*}
Then
     \begin{eqnarray}
     \log W_{n}^{n}\Rightarrow G  \quad \mbox{ in }\ell^{\infty}(\Psi),
      \end{eqnarray}
where  $G$ stands for a Gaussian element in $l^{\infty}(\Psi)$, each $d$-dimensional marginal $(G^{\psi_{1}},\cdots,G^{\psi_{d}})$ is a  normal random vector with mean
$$
\vec{\mu}=-\frac 12\big(V^{\psi_{i},\psi_{i}}, \, 1\le i\le d\big)
$$
and covariance structure
$$ \Sigma=\big(V^{\psi_{i},\psi_{j}}\big)_{1\le i,j\le d}.$$
To conclude the Introduction, two more remarks are given .
\begin{remark}
Observe in the discrete time case the Hellinger process can be computed as follows.
 \begin{eqnarray}
 h\Big(\frac 12;  P,P^{n,\psi}\Big)=2\sum_{i=1}^{n}E\big[(\xi^{n,\psi}_{i})^{2}|\mathcal{F}_{i-1}\big]
  \end{eqnarray}
    and
 \begin{eqnarray}
\imath(f_{1+\varepsilon},P, P^{n,\psi})=\sum_{i=1}^{n}E\big[(\xi_{i}^{n,\psi})^{2}1_{\{\frac{1}{1+\varepsilon}\le (\xi_{i}^{n,\psi}+1)^{2}<1+\varepsilon\}^{c}}\big|\mathcal{F}_{i-1}\big]
 \end{eqnarray}
Thus there is to some extent a similarity between our assumptions in Theorem \ref{m1}  and Nishiyama's assumptions.  However, it seems neater to use the Hellinger processes in continuous time case.

The integrability  condition ({\bf{Assumption \ref{a3}}}) of  partitioning entropy plays an important role in the proof of Theorem \ref{m1}. It is possible to  use the metric entropy condition, but we need to introduce a suitable pseudo-metric in the index set $\Psi$. The Hellinger processes would also be  very likely a good  candidate.
\end{remark}

\begin{remark}
It is rather interesting to consider the limiting behavior of the process $\log Z^{n}$ in $\ell^{\infty}([0,1]\times \Psi)$. To this end, we need to establish a  tightness criterion in the space $[0,1]\times \Psi$.  This is more complicate, and   will be left to the future work.

\end{remark}

\section{Proofs}

Let us start with a  decomposition. Observe that
\begin{eqnarray}
     N_t^{n,\psi} = \beta^{n,\psi} \cdot X_t^c +\Big(Y^{n,\psi}-1+\frac{\hat{Y}^{n,\psi}-a}{1-a}\Big)\ast (\mu-\nu)_t. \nonumber
\end{eqnarray}
It is easy to see
\begin{eqnarray}
     \Delta N_t^{n,\psi}=  \big(Y^{n,\psi}(t, \Delta X_t)-1\big)\mathbf{1}_{(\Delta X_t\neq 0)} -\frac{\hat{Y}^{n,\psi}_t-a_t}{1-a_t} \mathbf{1}_{(\Delta X_t=0)},\nonumber
\end{eqnarray}
and so we have
\begin{eqnarray*}
     \Delta L_t^{n,\psi} &=& \log \big(1+\Delta N_t^{n,\psi}\big) \\
     &=& \log \big(Y^{n,\psi}(t, \Delta X_t)\big)\mathbf{1}_{(\Delta X_t\neq 0)}+ \log\Big(1-\frac{\hat{Y}^{n,\psi}_t-a_t}{1-a_t} \Big)\mathbf{1}_{( \Delta X_t= 0)}.\nonumber
\end{eqnarray*}
Let $\mu^{n,\psi}$ be the jump measure of $L^{n,\psi}$ defined by
\begin{eqnarray*}
    \mu^{n,\psi}_t=   \sum_{s\le t}  \mathbf{1}_{(\Delta L_s^{n,\psi} \neq 0)} \varepsilon_{(s, \Delta L_s^{n,\psi})}
\end{eqnarray*}
and  $\nu^{n,\psi}$   the corresponding  predictable compensator. Then for any predictable function $W(\omega, t, x)$,
\begin{eqnarray*}
    W\ast \mu^{n,\psi}&=& \sum_{s\le t} W(\log Y^{n,\psi}(t, \Delta X_t))\mathbf{1}_{( \Delta X_t\neq 0)}\\
     & & +\sum_{s\le t} W\Big(\log\big(1-\frac{\hat{Y}^{n,\psi}-a}{1-a}\big)\Big) \mathbf{1}_{( \Delta X_t= 0)}\\
    &=& W(\log Y^{n,\psi})\ast \mu+ \sum_{s\le t} W\Big(\log\big(1-\frac{\hat{Y}^{n,\psi}_s-a_s}{1-a_s}\big)\Big) \mathbf{1}_{( \Delta X_s= 0)},
\end{eqnarray*}
and so by the fact that  $1-a_s$ is the predictable projection of $\mathbf{1}_{( \Delta X_s= 0)}$,
\begin{eqnarray}
   W\ast \nu^{n,\psi} = W(\log Y^{n,\psi})  \ast  \nu + \sum_{s\le t} W\Big(\log\big(1-\frac{\hat{Y}^{n,\psi}_s-a_s}{1-a_s}\big)\Big)  (1-a_s).
\end{eqnarray}
Given a positive number $\tau$, consider the truncation function
 \begin{equation}
       h_{\tau}(x)=x\mathbf{1}_{(|x|\le \tau)}
\end{equation}
and define
\begin{eqnarray}
      {\check{L}}_t^{n,\tau, \psi} = \sum_{s\le t} \big(\Delta L_s^{n,\psi}- h_\tau\big(\Delta L_s^{n,\psi} \big)  \big)
\end{eqnarray}
Thus combined together, we easily have a canonical decomposition
\begin{eqnarray}
       L_t^{n,\psi} &=&   {\check{L}}_t^{n,\tau, \psi} + \beta^{n,\psi}\cdot X_t^c-\frac 12 (\beta^{n,\psi})^2\cdot C_t  + h_\tau\ast (\mu^{n,\psi} -\nu^{n,\psi})_t \nonumber \\
         & & + \sum_{s\le t}\frac{\hat{Y}^{n,\psi}_s-a_s}{1-a_s}-\Big(Y^{n,\psi}-1+\frac{\hat{Y}^{n,\psi}-a}{1-a}\Big)\ast \nu_t \nonumber \\
         & &+ h_\tau(\log Y^{n,\psi})\ast \nu_t + \sum_{s\le t}\log\Big(1-\frac{\hat{Y}^{n,\psi}_s-a_s}{1-a_s}\Big) (1-a_s).\nonumber
\end{eqnarray}
For simplicity of writing, let
\begin{eqnarray}
        A^{n,\psi}_{1, t}= \beta^{n,\psi}\cdot X_{t}^c -\frac{1}{2}(\beta^{n,\psi})^2\cdot C_t,\nonumber
\end{eqnarray}
\begin{eqnarray}
        A^{n,\psi}_{2, t}=  h_\tau\ast (\mu^{n,\psi}-\nu^{n,\psi})_t,\nonumber
\end{eqnarray}
\begin{eqnarray}
        A^{n,\psi}_{3, t}&=&\sum_{s\le t}\frac{\hat{Y}^{n,\psi}_s-a_s}{1-a_s}-\Big(\hat{Y}^{n,\psi}-1+\frac{\hat{Y}^{n,\psi}-a}{1-a}\Big)\ast \nu_t \nonumber\\
         & &+ h_\tau\big(\log Y^{n,\psi}\big)\ast \nu_t + \sum_{s\le t}\log\Big(1-\frac{\hat{Y}^{n,\psi}_s-a_s}{1-a_s}\Big) (1-a_s),
\end{eqnarray}
Thus  we have
$$
 L_t^{n,\psi} =  {\check{L}}_t^{n,\tau, \psi} + A^{n,\psi}_{1, t} +A^{n,\psi}_{2, t}+A^{n,\psi}_{3, t}
$$
The proof of Theorem \ref{m1} will consist of a series of lemmas and propositions.
\begin{lemma} \label{le2}Under  Assumptions \ref{a1}, \ref{a2} and \ref{a3}, we have for each $\psi\in \Psi$ and $\tau>0$, as $n\rightarrow\infty$
\begin{eqnarray}\label{2.9}
\check{L}^{n,\tau, \psi}_t \stackrel P\longrightarrow 0
\end{eqnarray}
\end{lemma}
\begin{proof}
Set
\begin{eqnarray}
 \Upsilon^{n,\psi}_s=\frac{1-\hat{Y}^{n,\psi}_s}{1-a_s}. \label{3.55}
\end{eqnarray}
Then
\begin{eqnarray*}
 \check{L}^{n,\tau, \psi}_t = \log(Y^{n,\psi})\mathbf{1}_{(|\log(Y^{n,\psi})|>\tau)}\ast \mu_t +\sum_{s\le t}\log (\Upsilon^{n,\psi}_s )\mathbf{1}_{(|\log(\Upsilon^{n,\psi}_s)|>\tau)}\mathbf{1}_{(\Delta X_s=0)}.
\end{eqnarray*}
Also, for any $\delta>0$
\begin{eqnarray*}
& &  |\log(Y^{n,\psi})|\mathbf{1}_{(|\log(Y^{n,\psi})|>\tau)}\ast \nu_t \nonumber\\
  & &\le    |\log(Y^{n,\psi})|\mathbf{1}_{(|\log(Y^{n,\psi})|>\tau)}\mathbf{1}_{(|Y^{n,\psi}-1|\le \delta)}\ast \mu_t \nonumber\\
  & & \quad +  |\log(Y^{n,\psi})|\mathbf{1}_{(|\log(Y^{n,\psi})|>\tau)}\mathbf{1}_{(|Y^{n,\psi}-1|>\delta)}\ast \mu_t.
\end{eqnarray*}
By the Lenglart domination property (see page 35 of Jacod and Shiryaev \cite{js}),
\begin{eqnarray*}
 & & P\big(|\log(Y^{n,\psi})|\mathbf{1}_{(|\log(Y^{n,\psi})|>\tau)}\mathbf{1}_{(|Y^{n,\psi}-1|\le\delta)}\ast \mu_t> \varepsilon\big)  \nonumber\\
  & & \le  \frac{\eta}{\varepsilon}
    + P\big(|\log(Y^{n,\psi})|\mathbf{1}_{(|\log(Y^{n,\psi})|>\tau)}  \mathbf{1}_{(|Y^{n,\psi}-1|\le \delta)}\ast \nu_t> \eta\big).
\end{eqnarray*}
 Note for   $\delta<1$ there is a positive constant $c_\delta$ such that for any $x>0$
\begin{eqnarray}
|\log x| \mathbf{1}_{(|x-1|\le \delta)} \le c_\delta |x-1|\mathbf{1}_{(|x-1|\le \delta)},
\end{eqnarray}
so
\begin{eqnarray*}
   & &|\log(Y^{n,\psi})|\mathbf{1}_{(|\log(Y^{n,\psi})|>\tau)}\mathbf{1}_{(|Y^{n,\psi}-1|\le \delta)}\ast \nu_t \nonumber\\
  & & \le c_\delta | Y^{n,\psi}-1|\mathbf{1}_{(  e^{-\tau}<Y^{n,\psi}< e^\tau )^c}\mathbf{1}_{(|Y^{n,\psi}-1|\le \delta)}\ast \nu_t  \nonumber\\
  & & \stackrel P\longrightarrow 0.
\end{eqnarray*}
On the other hand, for each $\varepsilon\in (0,1)$
\begin{eqnarray*}
 & & P( |\log(Y^{n,\psi})|\mathbf{1}_{(|\log(Y^{n,\psi})|>\tau)}\mathbf{1}_{(|Y^{n,\psi}-1|> \delta)}\ast \mu_t >\varepsilon ) \nonumber\\
  & & \le P( \mathbf{1}_{(|Y^{n,\psi}-1|> \delta)}\ast \mu_t >\varepsilon ).
\end{eqnarray*}
Again, by the Lenglart domination property, it follows for any $\eta>0$
\begin{eqnarray*}
P( \mathbf{1}_{(|Y^{n,\psi}-1|> \delta)}\ast \mu_t >\varepsilon )&\le & \frac {\eta}{\varepsilon}+P( \mathbf{1}_{(|Y^{n,\psi}-1|> \delta)}\ast \nu_t >\eta )\nonumber\\
  & \le  & \frac {\eta}{\varepsilon}+P(|Y^{n,\psi}-1| \mathbf{1}_{(|Y^{n,\psi}-1|> \delta)}\ast \nu_t >\delta \eta ).
\end{eqnarray*}
Letting $n\rightarrow\infty$ and then $\eta\rightarrow 0$, we have
\begin{eqnarray*}
P( \mathbf{1}_{(|Y^{n,\psi}-1|> \delta)}\ast \mu_t >\varepsilon) \rightarrow 0.
\end{eqnarray*}
In combination, we have proved the desired statement.
\end{proof}

\begin{lemma} \label{le1}Under  Assumptions \ref{a1}, \ref{a2} and \ref{a3}, we have for each $\psi\in \Psi$, as $n\rightarrow\infty$

\begin{eqnarray}\label{2.5} (Y^{n,\psi}-1)\ast \nu_t  \stackrel P\longrightarrow 0,
\end{eqnarray}
\begin{eqnarray}\label{2.6}
(1-\sqrt {Y^{n,\psi}} )^2\ast \nu_t \stackrel P\longrightarrow 0,
\end{eqnarray}
\begin{eqnarray}\label{2.7} \sum_{s\le t}\Big(\sqrt{1-a_s}-\sqrt {1-\hat{Y}^{n,\psi}_s}\Big)^2 \stackrel P\longrightarrow 0.
\end{eqnarray}
Consequently,
\begin{eqnarray}\label{2.8}  \frac{1}{8}(\beta^{n,\psi})^2\cdot C_t \stackrel P\longrightarrow  V_t^{\psi}.
\end{eqnarray}
\end{lemma}
\begin{proof}
Obviously, for any $\varepsilon>0$
\begin{eqnarray*}
   | {Y}^{n,\psi}-1|\le \varepsilon+ | {Y}^{n,\psi}-1|\mathbf{1}_{(| {Y}^{n,\psi}-1|>\varepsilon)}.
\end{eqnarray*}
Also, by Assumption 2
\begin{eqnarray*}
    | {Y}^{n,\psi}-1|\mathbf{1}_{(| {Y}^{n,\psi}-1|>\varepsilon)}\ast \nu_t&\le&  \imath (h_{1+\varepsilon}; P, P^{n,\psi})\nonumber\\
    &\stackrel P\longrightarrow & 0 ,\quad n \rightarrow\infty.
\end{eqnarray*}
The desired (\ref{2.5}) holds.

Observe an elementary inequality: for any $0<\varepsilon<1$, there is a positive constant $c_\varepsilon$ such that
\begin{eqnarray*}
  (\sqrt{1+x}-1)^2\le \left\{
  \begin{array}{ll}
  |x|^2,  &  |x|\le \varepsilon<1,\\
  \\
  c_\varepsilon |x|,& |x|>\varepsilon.
  \end{array}
  \right.
\end{eqnarray*}
Then it follows
\begin{eqnarray*}
   (\sqrt{Y^{n,\psi}}-1)^2 \ast \nu_t&\le&  |Y^{n,\psi}-1|^2\mathbf{1}_{(|Y^{n,\psi}-1|\le \varepsilon)}\ast \nu_t\nonumber\\
   & &+ c_\varepsilon |Y^{n,\psi}-1|\mathbf{1}_{(|Y^{n,\psi}1|> \varepsilon)}\ast \nu_t,
\end{eqnarray*}
\begin{eqnarray*}
    |Y^{n,\psi}-1|^2\mathbf{1}_{(|Y^{n,\psi}-1|\le \varepsilon)}\ast \nu_t \le \varepsilon^2 \nu([0, t]\times \mathbb{R}),
\end{eqnarray*}
\begin{eqnarray*}
   |Y^{n,\psi}-1|\mathbf{1}_{(|Y^{n,\psi}-1|> \varepsilon)}\ast \nu_t \le \imath (h_{1+\varepsilon}; P, P^{n,\psi}).
\end{eqnarray*}
Thus under the Assumption \ref{a2},  we have by letting $n\rightarrow\infty$ and then   $\varepsilon\rightarrow 0$
\begin{eqnarray*}
   (\sqrt{Y^{n,\psi}}-1)^2 \ast \nu_t \stackrel P\longrightarrow 0.
\end{eqnarray*}

For (\ref{2.7}),   note
\begin{eqnarray*}
 \sum_{s\le t}\Big(\sqrt{1-a_s}-\sqrt {1-\hat{Y}^{n,\psi}_s}\Big)^2 &=& \sum_{s\le t}(1-a_s)\Big(1-\sqrt {\frac{1-\hat{Y}^{n,\psi}_s}{1-a_s}}\Big)^2\nonumber\\
 &=& \sum_{s\le t}(1-a_s)\Big(1-\sqrt {\Upsilon^{n,\psi}_s} \Big)^2,
\end{eqnarray*}
where $\Upsilon^{n,\psi}_s$ is as in (\ref{3.55}).
Then it easily follows
\begin{eqnarray*}
  \Big(1-\sqrt{\Upsilon^{n,\psi}_s}\Big)^2 \le |\Upsilon^\psi_s-1|^2\mathbf{1}_{(  |\Upsilon^{n,\psi}_s-1|\le \varepsilon)} +c_\varepsilon |\Upsilon^{n,\psi}_s-1| \mathbf{1}_{(  |\Upsilon^{n,\psi}_s-1|> \varepsilon)},
\end{eqnarray*}
\begin{eqnarray*}
  \sum_{s\le t}(1-a_s) |\Upsilon^{n,\psi}_s-1|^2\mathbf{1}_{(  |\Upsilon^{n,\psi}_s-1|\le \varepsilon)} \le \varepsilon^2 \sum_{s\le t}(1-a_s),
\end{eqnarray*}
\begin{eqnarray*}
  \sum_{s\le t}(1-a_s)   |\Upsilon^{n,\psi}_s-1| \mathbf{1}_{(  |\Upsilon^{n,\psi}_s-1|> \varepsilon)} \le  \imath (h_{1+\varepsilon}; P, P^{n,\psi}).
\end{eqnarray*}
Thus under the Assumption \ref{a2},  we have by first letting $n\rightarrow\infty$ and then   $\varepsilon\rightarrow 0$
\begin{eqnarray*}
    \sum_{s\le t}\Big(\sqrt{1-a_s}-\sqrt {1-\hat{Y}^{n,\psi}_s}\Big)^2 \stackrel P\longrightarrow 0.
\end{eqnarray*}
The proof is now complete.
\end{proof}

\begin{lemma} \label{le3}Under  Assumptions \ref{a1}, \ref{a2} and \ref{a3}, we have for each $\psi\in \Psi$ and $\tau>0$, as $n\rightarrow\infty$
\begin{eqnarray*}
 h_\tau\ast (\mu^{n,\psi}- \nu^{n,\psi})_t \stackrel P\longrightarrow 0.
\end{eqnarray*}
\end{lemma}

\begin{proof}
First, observe the quadratic variation of $h_\tau\ast (\mu^{n,\psi}- \nu^{n,\psi})_t$ is given by
\begin{eqnarray*}
 \langle h_\tau\ast (\mu^{n,\psi} - \nu^{n,\psi})\rangle_t &=& (h_\tau(\log Y^{n,\psi}))^2 \ast \nu_t \nonumber\\
 & &+ \sum_{s\le t} \Big(h_\tau \big( \log (1-\frac{\hat{Y}^{n,\psi}_s-a_s}{1-a_s}  ) \big)  \Big)^2(1-a_s).
\end{eqnarray*}
We shall prove that $\langle h_\tau\ast (\mu^{n,\psi} - \nu^{n,\psi})\rangle_t$ converges in probability to zero below. Note   there is a  $  \delta>0 $ such that for any $\varepsilon<\delta$
\begin{eqnarray*}
 \big|(h_\tau(\log(1+x)))^2 -4 (1-\sqrt{1+x})^2 \big|\le \left\{
 \begin{array}{ll}
  |x|^3, & |x|\le \varepsilon< \delta,\\
 c_\varepsilon |x|, & |x|>\varepsilon
 \end{array}
 \right.
\end{eqnarray*}
where $0<c_\varepsilon<\infty$. Thus it follows  for any $\varepsilon<\delta$
\begin{eqnarray*}
 (h_\tau(\log Y^{n,\psi}))^2 \ast \nu_t &\le &  | Y^\psi-1|^3\mathbf{1}_{(| Y^\psi-1|\le \varepsilon)} \ast \nu_t \nonumber\\
 & & +  c_\varepsilon | Y^{n,\psi}-1|\mathbf{1}_{(| Y^{n,\psi}-1|>\varepsilon)} \ast \nu_t.
\end{eqnarray*}
Hence letting    $n\rightarrow\infty$ and then $\varepsilon\rightarrow 0$ immediately yields
\begin{eqnarray*}
 (h_\tau(\log Y^{n,\psi}))^2 \ast \nu_t \stackrel P\longrightarrow 0.
\end{eqnarray*}
A similar argument shows
\begin{eqnarray*}
 \sum_{s\le t} \Big(h_\tau \big( \log (1-\frac{\hat{Y}^{n,\psi}_s-a_s}{1-a_s}  ) \big)  \Big)^2(1-a_s) \stackrel P\longrightarrow 0.
\end{eqnarray*}
Combined together, we have the desired statement.
\end{proof}

\begin{lemma}
 \label{le4}

Under  Assumptions \ref{a1}, \ref{a2} and \ref{a3}, we have for each $\psi\in \Psi$ and $\tau>0$, as $n\rightarrow\infty$
\begin{eqnarray*}
  A^{n,\psi}_{3, t} \stackrel P\longrightarrow 0.
\end{eqnarray*}
\end{lemma}
\begin{proof}
  Note there is a  $  \delta>0 $ such that for any $\varepsilon<\delta$
\begin{eqnarray*}
 \big| h_\tau(\log(1+x)) +2 (1-\sqrt{1+x})^2 -x\big|\le \left\{
 \begin{array}{ll}
  |x|^3, & |x|\le \varepsilon< \delta,\\
 c_\varepsilon |x|, & |x|>\varepsilon,
 \end{array}
 \right.
\end{eqnarray*}
where $0<c_\varepsilon<\infty$. Thus it follows  for any $\varepsilon<\delta$
\begin{eqnarray*}
 & & \Big|h_\tau\big(\log Y^{n,\psi}\big)\ast \nu_t+ 2(1-\sqrt{Y^{n,\psi}})^2\ast \nu_t -\Big(\hat{Y}^{n,\psi}-1+\frac{\hat{Y}^{n,\psi}-a}{1-a}\Big)\ast \nu_t \Big|\nonumber\\
 & &\le |Y^{n,\psi}-1|^3\mathbf{1}_{(| Y^{n,\psi}-1|\le \varepsilon)} \ast \nu_t + c_\varepsilon  |Y^{n,\psi}-1|\mathbf{1}_{(| Y^{n,\psi}-1|> \varepsilon)}\ast \nu_t.
\end{eqnarray*}
This implies
\begin{eqnarray*}
  \Big|h_\tau\big(\log Y^{n,\psi}\big)\ast \nu_t+ 2(1-\sqrt{Y^{n,\psi}})^2\ast \nu_t -\Big(\hat{Y}^{n,\psi}-1+\frac{\hat{Y}^{n,\psi}-a}{1-a}\Big)\ast \nu_t \Big|\stackrel P\longrightarrow 0.
\end{eqnarray*}
Since  it was proved $2(1-\sqrt{Y^{n,\psi}})^2\ast \nu_t\stackrel P\longrightarrow 0$, then we have
\begin{eqnarray*}
  \Big|h_\tau\big(\log Y^{n,\psi}\big)\ast \nu_t -\Big(\hat{Y}^{n,\psi}-1+\frac{\hat{Y}^{n,\psi}-a}{1-a}\Big)\ast \nu_t \Big|\stackrel P\longrightarrow 0.
\end{eqnarray*}
Similarly, we have
\begin{eqnarray*}
  \Big|\sum_{s\le t}\frac{\hat{Y}^{n,\psi}_s-a_s}{1-a_s} + \sum_{s\le t}h_\tau\Big(\log\big(1-\frac{\hat{Y}^{n,\psi}_s-a_s}{1-a_s}\big)\Big) (1-a_s) \Big|\stackrel P\longrightarrow 0.
\end{eqnarray*}
Combined together, the proof is complete.
\end{proof}
\begin{proposition}\label{np1}
  Under  Assumptions \ref{a1}, \ref{a2} and \ref{a3}, every finite-dimensional marginal of $L^{n}$ converges weakly.
\end{proposition}

\begin{proof}
From Lemmas \ref{le2}-\ref{le4},
         \begin{equation*}
            \beta^{n,\psi}_{s}\cdot X_{t}^{c}-\frac{1}{2} (\beta^{n,\psi}_{s})^{2}\cdot C_{t}
         \end{equation*}
has non-degenerate limiting finite-dimensional marginal laws, and the other part of $L^{n}$ asymptotically vanish.

For every $\psi$, the process
 \begin{equation*}
            \beta^{n,\psi}_{s}\cdot X_{t}^{c}-\frac{1}{2} (\beta^{n,\psi}_{s})^{2}\cdot C_{t}
 \end{equation*}
is a continuous semimartingale. Its predictable characterstics $(\hat{B}^{ \psi},\hat{C}^{ \psi}, 0)$ are
                  \begin{equation*}
                   \hat{B}^{n,\psi}_{t}=-\frac{1}{2} (\beta^{n,\psi} )^{2}\cdot C_{t},
                  \end{equation*}
                     \begin{equation*}
                   \hat{C}^{n,\psi}_{t}=(\beta^{n,\psi} )^{2}\cdot C_{t}.
                  \end{equation*}

By Lemma \ref{le1} and Assumption \ref{a2}, there is  a non-decreasing continuous function $V$, such that $V_{0}=0$,
  \begin{eqnarray}
  \sup_{\psi\in \Psi}|(\beta^{n,\psi} )^{2}\cdot C_{t} -V_t^{\psi,\psi}|\stackrel{P^{*}}{\rightarrow} 0,
  \end{eqnarray}
     \begin{eqnarray}
    (\beta^{n,\psi} \beta^{n,\psi} ) \cdot C_{t}\stackrel{P}{\rightarrow} V_t^{\psi,\phi}
    \end{eqnarray}
 for every $\psi,\phi\in \Psi$.

The proposition is now concluded by Theorem VIII.3.6 of Jacod and Shiryaev \cite{js}.
\end{proof}
Next we turn to verifying uniform tightness.
\begin{lemma} \label{le6}  Under  Assumptions \ref{a1}, \ref{a2} and \ref{a3}, we have for each   $\tau>0$, as $n\rightarrow\infty$
\begin{eqnarray}
       \sup_{\psi\in \Psi}\big|\check{L}_t^{n,\tau,\psi}\big|\stackrel P \longrightarrow 0.
\end{eqnarray}
\end{lemma}
\begin{proof}Recall
\begin{eqnarray}
       \check{L}_t^{n,\tau,\psi} &= & \sum_{s\le t} \log \big(Y^{n,\psi}(s, \Delta X_s)\big)\mathbf{1}_{(|\log Y^{n,\psi}(s, \Delta X_s|>\tau)}\mathbf{1}_{(\Delta X_s\neq 0)}\nonumber\\
       & & +\sum_{s\le t} \log \Big(1- \frac{\hat{Y}^{n,\psi}_s-a_s}{1-a_s}\Big)\mathbf{1}_{(|\log  (1- \frac{\hat{Y}^{n,\psi}_s-a_s}{1-a_s}  )|>\tau)}\mathbf{1}_{(\Delta X_s= 0)}.
\end{eqnarray}
Let us prove
\begin{eqnarray}
         \sup_{\psi\in \Psi}\sum_{s\le t} |\log \big(Y^{n,\psi}(s, \Delta X_s)\big)|\mathbf{1}_{(|\log Y^{n,\psi}(s, \Delta X_s)|>\tau)}\mathbf{1}_{(\Delta X_s\neq 0)} \stackrel P \longrightarrow 0 \label{Eq1}
\end{eqnarray}
and
\begin{eqnarray}
         \sup_{\psi\in \Psi}\sum_{s\le t} \log \Big(1- \frac{\hat{Y}^{n,\psi}_s-a_s}{1-a_s}\Big)\mathbf{1}_{(|\log  (1- \frac{\hat{Y}^{n,\psi}_s-a_s}{1-a_s}  )|>\tau)}\mathbf{1}_{(\Delta X_s= 0)} \stackrel P \longrightarrow 0. \label{Eq2}
\end{eqnarray}
For (\ref{Eq1}), note
\begin{eqnarray}
     & &    |\log \big(Y^{n,\psi}(s, \Delta X_s)\big)|\mathbf{1}_{(|\log Y^{n,\psi}(s, \Delta X_s)|>\tau)}\nonumber\\ & & =\log \big(Y^{n,\psi}(s, \Delta X_s)\big) \mathbf{1}_{( \log Y^{n,\psi}(s, \Delta X_s) >\tau)}\nonumber\\
         &   &\quad -\log \big(Y^{n,\psi}(s, \Delta X_s)\big) \mathbf{1}_{( \log Y^{n,\psi}(s, \Delta X_s) <-\tau)}.\nonumber
\end{eqnarray}
Thus we need only to prove
\begin{eqnarray}
      \sum_{s\le t} \sup_{\psi\in \Psi} \log \big(Y^{n,\psi}(s, \Delta X_s)\big) \mathbf{1}_{( \log Y^{n,\psi}(s, \Delta X_s) >\tau)}\mathbf{1}_{(\Delta X_s\neq 0)}  \stackrel P \longrightarrow 0 \label{Eq3}
\end{eqnarray}
and
\begin{eqnarray}
     \sum_{s\le t} \inf_{\psi\in \Psi} \log \big(Y^{n,\psi}(s, \Delta X_s)\big) \mathbf{1}_{( \log Y^{n,\psi}(s, \Delta X_s) <-\tau)}\mathbf{1}_{(\Delta X_s\neq 0)}  \stackrel P \longrightarrow 0. \label{Eq4}
\end{eqnarray}
Let us first look at (\ref{Eq3}). Set
\begin{eqnarray}
      \overline{\mathbf{Y}}^{n,\Psi}(\omega; s, x)= \Big[\sup_{\psi\in \Psi} Y^{n,\psi}(\omega; s, x)\Big]_{\tilde{\mathcal{P}}, M_{\mu}^{P}}.
\end{eqnarray}
Then
\begin{eqnarray}
  & &  \sum_{s\le t}   \sup_{\psi\in \Psi} \log \big(Y^{n,\psi}(s, \Delta X_s)\big) \mathbf{1}_{( \log Y^{n,\psi}(s, \Delta X_s) >\tau)}\mathbf{1}_{(\Delta X_s\neq 0)}\nonumber\\
  & &\le  \sum_{s\le t}    \log \big(\overline{\mathbf{Y}}^{n,\Psi}(\omega; s, \Delta X_s)\big) \mathbf{1}_{(\log \overline{\mathbf{Y}}^{n,\Psi}(\omega; s, \Delta X_s) >\tau)}\mathbf{1}_{(\Delta X_s\neq 0)}\nonumber\\
  & &=   \log \big(\overline{\mathbf{Y}}^{n,\Psi} \big) \mathbf{1}_{(\log \overline{\mathbf{Y}}^\Psi  >\tau)} \ast \mu_t.
\end{eqnarray}
For any $\varepsilon, \eta>0$,
\begin{eqnarray}
& &P(\log \big(\overline{\mathbf{Y}}^{n,\Psi} \big) \mathbf{1}_{(\log \overline{\mathbf{Y}}^{n,\Psi}  >\tau)} \ast \mu_t>\varepsilon)\nonumber\\
& & \le  \frac{\eta}{\varepsilon}+ P(\log \big(\overline{\mathbf{Y}}^\Psi \big) \mathbf{1}_{(\log \overline{\mathbf{Y}}^{n,\Psi}  >\tau)} \ast \nu_t>\eta).\nonumber
\end{eqnarray}
Note  $x-1> \log x$ if $\log x>0$. Then
\begin{eqnarray}
 \log \big(\overline{\mathbf{Y}}^{n,\Psi} \big) \mathbf{1}_{(\log \overline{\mathbf{Y}}^{n,\Psi}  >\tau)}\le (\overline{\mathbf{Y}}^{n,\Psi}-1)\mathbf{1}_{( \overline{\mathbf{Y}}^{n,\Psi}-1  >\tau)}.
\end{eqnarray}
Recall the definition of $\imath(f_{1+\varepsilon},P,P^{n,\psi})$ and Assumptions \ref{a1} and \ref{a2}, we can obtain (\ref{Eq3}). The proofs of $(\ref{Eq4})$ and $(\ref{Eq2})$  are similar.
 \end{proof}

\begin{lemma}\label{le7} Under  Assumptions \ref{a1}, \ref{a2} and \ref{a3}, for any $\varepsilon, \eta>0$, there is a $\delta>0$ and a partition $\Pi(\delta)=\{\Psi(\delta),\, 1\le k\le N(\delta)\}$ such that
\begin{eqnarray}
  \limsup_{n\rightarrow\infty}  P\Big(\sup_{1\le k\le N_{\Pi }(\delta)}\sup_{\phi, \psi\in \Psi(\delta, k)}|A_t^{n,2, \psi}- A_t^{n,2, \phi} |>\varepsilon\Big)\le \eta.
\end{eqnarray}
\end{lemma}

\begin{proof}   Let us fix $\varepsilon, \eta>0$. First, note
\begin{eqnarray*}
  A_t^{n,2, \psi} &=&  h_\tau(\log Y^{n,\psi})\ast (\mu -\nu )_t   \\
  & & + \sum_{s\le t} h_\tau(\log (1- \frac{\hat{Y_s}^{n,\psi}-a_s}{1-a_s}))\mathbf{1}_{(\Delta X_s=0)} \\
  & &-h_\tau(\log (1- \frac{\hat{Y_s}^{n,\psi}-a_s}{1-a_s}))(1-a_s).
\end{eqnarray*}
Let
\begin{eqnarray*}
   J^{n,1, \psi}_t= h_\tau(\log Y^{n,\psi})\ast (\mu -\nu )_t
\end{eqnarray*}
and
\begin{eqnarray*}
J^{n, 2, \psi}_t &=& \sum_{s\le t}\Big[ h_\tau\Big(\log (1-\frac{\hat{Y}^{n,\psi}_s-a_s}{1-a_s})\Big)\mathbf{1}_{(\Delta X_s=0)} \nonumber \\
& &\qquad - h_\tau\Big(\log (1-\frac{\hat{Y}^{n,\psi}_s-a_s}{1-a_s})\Big)(1-a_s) \Big].
\end{eqnarray*}
We shall treat $J^{n,1, \psi}_t$ and $J^{n,2, \psi}_t$ separately below. Let us only focus on the $J^{n,2, \psi}_t$ since the $J^{n, 1, \psi}_t$ is similar and simpler.

According to Assumption \ref{a3}, there is a sufficiently large positive finite constant $K$ such that
\begin{eqnarray}\label{nr3}
  \limsup_{n\rightarrow\infty} P(\|\mathfrak{H}^n\|_{\Pi, t}>K)\le \frac{\eta}{4}.
\end{eqnarray}
Thus we only need to condition on the event $\{\mathfrak{H}^n\|_{\Pi, t}>K\}$.   In particular, we shall  prove
     \begin{equation}\label{r1}
      E^{*} \max_{1\le k\le N_{\Pi}(\delta)}\sup_{\psi,\phi\in \Psi(\delta,k)}|J_{t}^{2, \psi }-J_{t}^{2, \phi} |\mathbf{1}_{\{||\mathfrak{H}||_{\Pi,t}\le K\}}\le c_{11}\int_{0}^{\delta}H_{\Pi}(\varepsilon)d\varepsilon.
 \end{equation}
Assuming (\ref{r1}), we can take $\delta$ so  small that
 \begin{equation}\label{r2}
        c_{11}\int_{0}^{\delta}H_{\Pi}(\varepsilon)d\varepsilon< \frac{\varepsilon \eta}{4},
 \end{equation}
 from which it follows by the Markov inequality
 \begin{eqnarray*}
 \limsup_{n\rightarrow\infty}  P(\sup_{1\le k\le N_{\Pi }(\delta)}\sup_{\phi, \psi\in \Psi(\delta, k)}|J_t^{n,2, \psi}- J_t^{n,2, \phi} |>\varepsilon, ||\mathfrak{H}||_{\Pi,t}\le K )\le \frac{\eta}{4}
 \end{eqnarray*}
 This in turn  together with (\ref{nr3}) implies
 \begin{eqnarray*}
 \limsup_{n\rightarrow\infty}  P\Big(\sup_{1\le k\le N_{\Pi }(\delta)}\sup_{\phi, \psi\in \Psi(\delta, k)}|J_t^{n,2, \psi}- J_t^{n,2, \phi} |>\varepsilon \Big)\le \frac{\eta}{2}.
 \end{eqnarray*}
It remains to prove (\ref{r1}). For every integer $p\ge 0$, construct a  nested refinement partition $\Pi(2^{-p}\delta)=\{\Psi(2^{-p}\delta; k), 1\le k\le N_\Pi(2^{-p}\delta)\} $ of $\Psi$, and then choose an element $\psi_{p,k}$ from each partitioning set $\Psi(2^{-p}\delta;k)$ in such a way that
\begin{eqnarray}
\{\psi_{p,k}:1\le k\le N_{\Pi}(2^{-p}\delta)\}\subset\{\psi_{p+1,k}:1\le k\le N_{\Pi}(2^{-p-1}\delta)\}.
\end{eqnarray}
For every $\psi\in \Psi$ and each $p\ge 0$,  define $\pi_{p} \psi=\psi_{p,k}$  and $\Pi_{p}\psi=\Psi(2^{-p}\delta;k)$ whenever $\psi\in \Psi(2^{-p}\delta;k)$. Obviously, $ \Pi_{p}\psi\subseteq\Pi_{p-1}\psi $. Define
\begin{eqnarray*}
 W(\Pi_{p}\psi)_t &=& \big[\sup_{\varphi,\phi\in \Pi_{p}\psi}|h_\tau(\log(1-(\frac{\hat{Y}_{s}^{n,\varphi}-a_{s}}{1-a_{s}})))
\nonumber\\
& & \quad-h_\tau(\log(1-(\frac{\hat{Y}_{s}^{n,\phi}-a_{s}}{1-a_{s}})))|\big]_{\tilde{\mathcal{P}}, M_{\mu}^{P}}.
\end{eqnarray*}
Note $W(\Pi_{p}\psi)\le W(\Pi_{p-1}\psi)$.  Set
 \begin{eqnarray}\label{ner1}
\alpha_{p}=\frac{2^{-p+1}\delta  }{H_{\Pi}(2^{-p-1}\delta)}K, \quad p\ge 0
\end{eqnarray}
and
\begin{eqnarray*}
A_{0}(\psi)=\mathbf{1}_{\{W(\Pi_{0}\psi)\le \alpha_{0}\}}, \quad B_{0}(\psi)=\mathbf{1}_{\{W(\Pi_{0}\psi)> \alpha_{0}\}}
\end{eqnarray*}
and for $p\ge 1$
\begin{eqnarray*}
A_{p}(\psi)=\mathbf{1}_{\{W(\Pi_{0}\psi)\le \alpha_{0},\cdots,W(\Pi_{p-1}\psi)\le \alpha_{p-1},W(\Pi_{p}\psi)\le \alpha_{p}\}},
\end{eqnarray*}
\begin{eqnarray*}
B_{p}(\psi)=\mathbf{1}_{\{W(\Pi_{0}\psi)\le \alpha_{0},\cdots,W(\Pi_{p-1}\psi)\le \alpha_{p-1},W(\Pi_{p}\psi)> \alpha_{p}\}}.
\end{eqnarray*}
It is easy to see
\begin{eqnarray*}
 A_0+B_0=1
\end{eqnarray*}
and  for each $p\ge 1$
\begin{eqnarray*}
 A_p+B_p=A_{p-1}.
\end{eqnarray*}
Hence it follows for any $q\ge 1$
\begin{eqnarray*}
  A_q+B_q+B_{q-1}+\cdots+B_0=1.
\end{eqnarray*}

Note $\pi_{0}\psi=\pi_{0}\phi$ if $ \psi, \phi\in \Psi(\delta, k)  $, and so we have
\begin{eqnarray*}
\sup_{\phi, \psi\in \Psi(\delta, k)}|J_{t}^{n,2, \psi }-J_{t}^{n,2,  \phi }|\le 2\sup_{\psi\in \Psi} |J_{t}^{n,2, \psi }-J_{t}^{n,2, \pi_{0}\psi }|.
\end{eqnarray*}
 It is now enough to show
\begin{equation}\label{r2}
       E^{*} \sup_{\psi\in \Psi}|J_{t}^{n,2, \psi }-J_{t}^{n,2, \pi_{0}\psi }|1_{\{||\mathfrak{H}||_{\Pi,t}\le K\}}\le c_{12}\int_{0}^{\delta}H_{\Pi}(\varepsilon)d\varepsilon.
\end{equation}
We have the following identity
\begin{eqnarray}
  J_{t}^{n,2, \psi }-J_{t}^{n,2, \pi_{0}\psi }&=&  \sum_{p=0}^q (J_{t}^{n,2, \psi }-J_{t}^{n,2, \pi_{p}\psi })B_p(\psi)\label{J4}\\
  & & +  (J_{t}^{n,2, \psi }-J_{t}^{n,2, \pi_{q}\psi })A_q(\psi)\nonumber\\
  & & + \sum_{p=1}^q (J_{t}^{n,2, \pi_p\psi }-J_{t}^{n,2, \pi_{p-1}\psi })A_{p-1}(\psi).\nonumber
\end{eqnarray}
Denote
     \begin{eqnarray*}
     M^{n,p,\psi}_{t}=(J_{t}^{n,2,\psi }-J_{t}^{n,2, \pi_{p}\psi})B_{p}(\psi).
     \end{eqnarray*}
 We shall only establish (\ref{r2}) for $M^{n,p,\psi}$ since the other three terms in RHS of  (\ref{J4}) can be similarly treated.

 Obviously, $M^{n,p,\psi}$ is a local martingale, and
 \begin{eqnarray*}
 |\Delta M^{n,p,\psi}_t|\le |W(\Pi_{p-1}(\psi))_t|B_{p}(\psi)\le \alpha_{p-1}.
 \end{eqnarray*}
On the other hand, the predictable quadratic variation of $M^{n,p,\psi}_t$ satisfies
\begin{eqnarray*}
  \langle M^{n,p,\psi}  \rangle_t = \sum_{s\le t} \Big[h_\tau(\log (1-\frac{\hat{Y}^{n,\psi}_s-a_s}{1-a_s}))-h_\tau(\log (1-\frac{\hat{Y}^{n,\pi_{p}\psi}_s-a_s}{1-a_s}))\Big]^2(1-a_s).
 \end{eqnarray*}
By an elementary calculation, for any $\gamma\in (0,1)$ there is a constant $c_\gamma$ such that
  \begin{eqnarray*}
  |\log x-\log y|\le c_\gamma|\sqrt{x}-\sqrt{y}|
  \end{eqnarray*}
   whenever $x,y\in [1-\gamma,1+\gamma]$. Then it follows
       \begin{eqnarray*}
       \langle M^{n,p,\psi}\rangle_{1}&\le&  \sum_{s\le t} \Big[ \sqrt{1-\frac{\hat{Y}^{n,\psi}_s-a_s}{1-a_s}}- \sqrt{ (1-\frac{\hat{Y}^{n,\pi_{p}\psi}_s-a_s}{1-a_s})}\Big]^2(1-a_s)\nonumber\\
       &\le & \sum_{s\le t} \Big[ \sqrt{1- \hat{Y}^{n,\psi}_s }- \sqrt{  1- \hat{Y}^{n,\pi_{p}\psi}_s}\Big]^2 \nonumber\\
       &\le & 2 h(\frac 12;  P^{n,\psi}, P^{n,\pi_{p}\psi})\le 2(2^{-p }\delta)^2 ||\mathfrak{H}||^2_{\Pi,t}.
        \end{eqnarray*}
By Bernstein-Freedman's inequality (see  Lemma 3.2 of Nishiyama \cite{n1}) for local martingale with bounded jumps,  it follows for $\varepsilon>0$,
        \begin{eqnarray*}
      P( |M^{n,p,\psi}_{t}|>\varepsilon,||\mathfrak{H}||_{\Pi,t}\le K )\le 2\exp\Big(-\frac{\varepsilon^{2}}{2[\alpha_{p-1}
        \varepsilon+(2^{-p}\delta K)^{2}]}\Big).
        \end{eqnarray*}
This, in turn together with   Lemma 2.2.10 of van der Vaart  and Wellner \cite{vw}, yields
      \begin{eqnarray*}
      & &E^{*} \sup_{\psi\in \Psi}|M^{n,p,\psi}|\mathbf{1}_{\{||\mathfrak{H}||_{\Pi,1}\le K\}}\\
      & & \le \alpha_{p-1} H_{\Pi}(2^{-p}\delta) ^{2}+2^{-p}\delta c_{13} H_{\Pi}(2^{-p}\delta)\\
      & & \le  c_{13}2^{-p}\delta  H_{\Pi}(2^{-p}\delta)
      \end{eqnarray*}
and
      \begin{eqnarray*}
    & & E^{*} \sup_{\psi\in \Psi}  \Big|\sum_{p=0}^q  J_{t}^{n,2, \psi }-J_{t}^{n,2, \pi_{p}\psi }\Big|B_p(\psi)\mathbf{1}_{\{||\mathfrak{H}||_{\Pi,1}\le K\}}\\
    & & \le c_{13}\sum_{p=0}^{q}2^{-p}\delta  H_{\Pi}(2^{-p}\delta)\le  c_{14}\int_{0}^{\delta}H_{\Pi}(\varepsilon)d\varepsilon.
      \end{eqnarray*}
Thus (\ref{r1}) is obtained, and so complete the proof.
 \end{proof}

\begin{lemma}\label{le8} Under  Assumptions \ref{a1}, \ref{a2} and \ref{a3}, for any $\varepsilon, \eta>0$, there is a $\delta>0$ and a partition $\Pi(\delta)=\{\Psi(\delta),\, 1\le k\le N(\delta)\}$ such that
\begin{eqnarray*}
 \limsup_{n\rightarrow\infty}P\Big(\sup_{1\le k\le N_{\Pi }(\delta)}\sup_{\phi, \psi\in \Psi(\delta, k)}|A_t^{n,1, \psi}- A_t^{n,1, \phi} |>\varepsilon\Big)\le \eta.
\end{eqnarray*}
\end{lemma}

\begin{proof} Recall
\begin{eqnarray}
        A^{n,1, \psi}_{ t}= \beta^{n,\psi}\cdot X_{t}^c -\frac{1}{2}(\beta^{n,\psi})^2\cdot C_t.\nonumber
\end{eqnarray}
It is enough to prove the following two statements
\begin{eqnarray}
 \quad\quad\limsup_{n\rightarrow\infty}P\Big(\sup_{1\le k\le N_{\Pi }(\delta)}\sup_{\phi, \psi\in \Psi(\delta, k)}|(\beta^{n,\psi})^2\cdot C_t- (\beta^{n,\psi})^2\cdot C_t |>\varepsilon\Big)\le \frac {\eta}{2} \label{A1-1}
\end{eqnarray}
and
\begin{eqnarray}
 \limsup_{n\rightarrow\infty}P\Big(\sup_{1\le k\le N_{\Pi }(\delta)}\sup_{\phi, \psi\in \Psi(\delta, k)}|\beta^{n,\psi}\cdot X_{t}^c- \beta^{n,\psi}\cdot X_{t}^c |>\varepsilon\Big)\le \frac {\eta}{2}.\label{A1-2}
\end{eqnarray}
We shall concentrate on proving (\ref{A1-1}) below since (\ref{A1-2}) is similar. The proof is completely similar to that of Lemma \ref{le7} with some minor modifications.  For every integer $p\ge 0$,  choose an element $\psi_{p,k}$ from each partitioning set $\Psi(2^{-p}\delta;k)$ in such a way that
\begin{eqnarray*}
\{\psi_{p,k}:1\le k\le N_{\Pi}(2^{-p}\delta)\}\subset\{\psi_{p+1,k}:1\le k\le N_{\Pi}(2^{-p-1}\delta)\}.
\end{eqnarray*}
and  define $\pi_{p} \psi=\psi_{p,k}$  and $\Pi_{p}\psi=\Psi(2^{-p}\delta;k)$ whenever $\psi\in \Psi(2^{-p}\delta;k)$. Note
\begin{eqnarray*}
& &\sup_{1\le k\le N_{\Pi }(\delta)}\sup_{\phi, \psi\in \Psi(\delta, k)}|(\beta^{n,\psi})^2\cdot C_t- (\beta^{n,\phi})^2\cdot C_t | \nonumber\\
& &\le 2 \sup_{\psi\in \Psi}| (\beta_{s}^{n,\psi})^{2}\cdot C_t - (\beta_{s}^{n,\pi_{0}\psi})^{2} \cdot C_{t}|,
\end{eqnarray*}
then a main step is to prove
\begin{equation}\label{r3}
       E^{*} \sup_{\psi\in \Psi}| (\beta^{n,\psi})^{2}\cdot C_t - (\beta^{n,\pi_{0}\psi})^{2} \cdot C_{t}|\mathbf{1}_{\{||\Xi||_{\Pi,1}\le K\}}\le c_{14}\int_{0}^{\delta}H_{\Pi}(\varepsilon)d\varepsilon.
     \end{equation}
To this end, for $p \ge 0$, set
\begin{eqnarray*}
\Gamma(\Pi_{p}\psi)_t= \big[\sup_{\varphi,\phi\in \Pi_{p}\psi}|(\beta_{s}^{n,\varphi})^{2}-(\beta_{s}^{n,\phi})^{2}|\big]_{\tilde{\mathcal{P}},M_{\mu}^{P}}.
\end{eqnarray*}
Obviously, $\Gamma(\Pi_{p}\psi)\le \Gamma(\Pi_{p-1}\psi)$. Define
\begin{eqnarray*}
D_{0}(\psi)=1_{\{\Gamma(\Pi_{0}\psi)\le \alpha_{0}\}}, \quad E_{0}(\psi)=1_{\{\Gamma(\Pi_{0}\psi)> \alpha_{0}\}},
\end{eqnarray*}
\begin{eqnarray*}
D_{p}(\psi)=1_{\{\Gamma(\Pi_{0}\psi)\le \alpha_{0},\cdots,\Gamma(\Pi_{p-1}\psi)\le \alpha_{p-1},\Gamma(\Pi_{p}\psi)\le \alpha_{p}\}},
\end{eqnarray*}
\begin{eqnarray*}
E_{p}(\psi)=1_{\{\Gamma(\Pi_{0}\psi)\le \alpha_{0},\cdots,\Gamma(\Pi_{p-1}\psi)\le \alpha_{p-1},\Gamma(\Pi_{p}\psi)> \alpha_{p}\}},
\end{eqnarray*}
where $\alpha_{p}$ is as in (\ref{ner1}).

Note we have the following identity
 \begin{eqnarray*}
    (\beta^{n,\psi})^{2}\cdot C_t-(\beta^{n,\pi_{0}\psi})^{2} \cdot C_{t}
   &= & \sum_{p=1}^{q}  ((\beta^{n,\pi_{p}\psi})^{2} \cdot C_{t} -(\beta^{n,\pi_{p-1}\psi})^{2} \cdot C_{t} ) D_{p-1}(\psi)\\
   & &+ ((\beta^{n,\psi})^{2} \cdot C_{t} -(\beta^{n,\pi_{q}\psi})^{2} \cdot C_{t}) D_{q}(\psi)\\
      &  &+ \sum_{p=0}^{q} ((\beta^{n,\psi})^{2} \cdot C_{t} -(\beta^{n,\pi_{p}\psi})^{2} \cdot C_{t})  E_{p}(\psi).
     \end{eqnarray*}
and
       \begin{eqnarray*}
        & &  \big((\beta^{n,\psi})^{2} \cdot C_{t} -(\beta^{n,\pi_{p}\psi})^{2} \cdot C_{t}\big)E_{p}(\psi)\\
        &&\le   \Gamma(\Pi_{p}\psi) (C_{t}-C_{0})E_{p}(\psi)  \\
        & &\le \frac{\Gamma(\Pi_{p}\psi) ^2}{\alpha_p} (C_{t}-C_{0})E_{p}(\psi).
       \end{eqnarray*}
 In addition, it is easy to see
 \begin{eqnarray*}
        & &E^{*} \sup_{\psi\in \Psi}((\beta^{n,\pi_p \psi})^{2} \cdot C_{t}-(\beta^{n,\pi_{p-1} \psi})^{2} \cdot C_{t})D_{p}(\psi)\\
        & & \le 2^{-p}\delta H_{\Pi}(2^{-p}\delta)
       \end{eqnarray*}
 by  Schwarz's inequality.  Thus
 \begin{eqnarray*}
     & & E^{*} \sup_{\psi\in \Psi}| (\beta^{n,\psi})^{2}\cdot C_t - (\beta^{n,\pi_{0}\psi})^{2} \cdot C_{t}|\mathbf{1}_{\{||\mathfrak{H}||_{\Pi,1}\le K\}}\\
   &  & \le \sum_{p=1}^{q} K_{2} 2^{-p-1}\delta
   H_{\Pi}(2^{-p-1}\delta)\\
   & &\le  K_{2} \int_{0}^{\delta}H_{\Pi}(\varepsilon)d\varepsilon.
     \end{eqnarray*}

Thus, (\ref{A1-1}) is proved. We complete the proof of this lemma.
\end{proof}
We can obtain the following proposition by Lemmas \ref{le6} - \ref{le8}.
\begin{proposition}\label{np2}
  Under  Assumptions \ref{a1}, \ref{a2} and \ref{a3}, for any $\varepsilon, \eta>0$, there is a $\delta>0$ and a partition $\Pi(\delta)=\{\Psi(\delta),\, 1\le k\le N(\delta)\}$ such that
\begin{eqnarray*}
 \limsup_{n\rightarrow\infty}P\Big(\sup_{1\le k\le N_{\Pi }(\delta)}\sup_{\phi, \psi\in \Psi(\delta, k)}|L^{n,\psi}- L^{n,\phi} |>\varepsilon\Big)\le \eta.
\end{eqnarray*}
\end{proposition}

\noindent
{\bf The proof of Theorem \ref{m1}.} Proposition \ref{np2} implies the asymptotic equicontinuity of $L^n$, and the asymptotic marginal distribution of  $L^n$ is obtained by Proposition \ref{np1}. Then we can obtain Theorem \ref{m1} by these two propositions and Theorem 1.1 in Nishiyama \cite{n3}.

\section*{Acknowledgments}
The authors would like to thank the anonymous referees and the Associate Editor for careful reading and constructive comments.
This work  was supported by the National Natural Science Foundation of China (No.11371317, 11701331,  11731012, 11871425) ,  Fundamental Research Funds for Central Universities, Shandong Provincial Natural Science Foundation (No. ZR2017QA007) and Young Scholars Program of Shandong University.

\vskip3mm

\bibliographystyle{amsplain}

\end{document}